\documentclass[a4paper,fleqn]{cas-dc}
\usepackage[authoryear,longnamesfirst]{natbib}
\usepackage{amsmath}
\usepackage[ruled,linesnumbered]{algorithm2e}
\usepackage{algorithmic}
\usepackage{makecell}

\def\tsc#1{\csdef{#1}{\textsc{\lowercase{#1}}\xspace}}
\tsc{WGM}
\tsc{QE}

\begin{document}
	\let\WriteBookmarks\relax
	\def\floatpagepagefraction{1}
	\def\textpagefraction{.001}
	\shorttitle{}
	\shortauthors{Yajie Wen et~al.}
	
	\title [mode = title]{A Block-Based Heuristic Algorithm for the Three-Dimensional Nuclear Waste Packing Problem} 
	                     
\author[1]{Yajie Wen}[
orcid = 0009-0004-2053-0583,
]
\ead{wynne.jei@gmail.com}
\fnmark[1]

\author[2]{Defu Zhang}[
%orcid = 0009-0004-2053-0583,
]
\fnmark[2]
\ead{dfzhang@xmu.edu.cn}

\affiliation{organization={Xiamen University},
	addressline={Siming south street 443}, 
	city={Xiamen}, 
	postcode={361000}, 
	state={Fujian},
	country={China}}
	
	\begin{abstract}
In this study, we present a block-based heuristic search algorithm to address the nuclear waste container packing problem in the context of real-world nuclear power plants. Additionally, we provide a dataset comprising 1600 problem instances for future researchers to use. Experimental results on this dataset demonstrate that the proposed algorithm effectively enhances the disposal pool's space utilization while minimizing the radiation dose within the pool. The code and data employed in this study are publicly available to facilitate reproducibility and further investigation.
	\end{abstract}

%	\begin{highlights}
%		\item Research highlights item 1
%		\item Research highlights item 2
%		\item Research highlights item 3
%	\end{highlights}
	
	\maketitle
	
	\section{Introduction}
	This paper addresses the challenge of nuclear waste treatment in the context of a real-life nuclear power plant. During regular operation, nuclear power plants generate radioactive waste, which must be collected, classified into distinct boxes, and stored in pre-dug disposal pools. The waste packages are sealed with cement to ensure the safety of the environment and human health. 
	
	The waste packages within the pools must be strategically arranged to maximize disposal pool usage and reduce radiation exposure per unit time (measured in Sieverts per second) in each pool. This problem can be formulated as a three-dimensional knapsack problem with an additional constraint on the minimum radiation limit. To tackle this, the paper proposes a block-based heuristic algorithm, BSNA, to identify optimal placement strategies for the waste boxes. 
	
	Furthermore, we generate 1600 synthetic nuclear waste packing problems based on real-world nuclear waste data to evaluate the proposed algorithm's performance. These generated problems are intentionally more complex than typical real-world scenarios, incorporating larger disposal pools and a greater variety and volume of nuclear waste boxes. This provides a comprehensive test bed for assessing the algorithm's effectiveness in more challenging settings.
	\section{Problem Description}
	The dose rate of a nuclear waste box at a given position, measured in Sieverts per second, is determined by the number of radioactive decay occurring per second, known as the activity, multiplied by the dose rate constant of the radioactive element. This value is then divided by the square of the distance between the nuclear element and the observation point. Since the dose rate constant for a specific nuclear element is a fixed value, the radiation dose of the nuclear waste box is directly proportional to the activity and inversely proportional to the square of the distance. The relationship is mathematically expressed as (\ref{hc}).
	\begin{eqnarray}\label{hc}
		\dot{H} = \frac{\Gamma \cdot A}{r^2}
	\end{eqnarray}

\begin{itemize}
	\item Dose Rate \( \dot{H} \) \([\text{Sv}/t]\): The amount of radiation absorbed per unit time, adjusted for biological effects.
\item Dose Rate Constant \( \Gamma \) \([\text{Sv} \cdot \text{m}^2/(\text{Bq} \cdot t)]\): A factor that describes how much radiation dose is received per unit of radioactive source activity at a given distance.
\item Activity \( A \) \([\text{Bq}]\): The number of radioactive decay occurring per second.
\item Distance \( r \) \([\text{m}]\): The distance from the radiation source to the point of measurement.
\end{itemize}
	
	For simplicity, let the vertical distance from the center point of the nuclear waste box to the top of the disposal pool be \(r\). Given a disposal pool with length \(L\), width \(W\), and height \(H\), and a set of \(n\) nuclear waste boxes, each with length \(l_i\), width \(w_i\), height \(h_i\), activity \(A_i\), and dose rate constant \(\Gamma_i\), the objective of the problem is to place as many nuclear waste boxes as possible into the disposal pool to maximize the utilization rate of the disposal pool, while minimizing the sum of the dose rates of all the nuclear waste boxes in the disposal pool. This problem is similar to a variant of the three-dimensional bin packing problem, as it has an additional optimization objective.
	
	\section{Literature Review}
	
	The three-dimensional bin packing problem (3D-BPP or 3D-CLP) is a well-known NP-hard optimization problem. Over the years, various approaches have been developed to tackle this problem, including exact algorithms, heuristic algorithms, metaheuristic algorithms, and more recently, deep learning-based methods. Among these, heuristic and metaheuristic algorithms are the most widely adopted due to their practical efficiency in solving large-scale problems.
	
	\subsection{Exact Algorithms}
	Exact algorithms aim to obtain optimal solutions for 3D-BPP. However, due to the complexity and constraints of the problem, these algorithms are typically suitable for small-scale instances. Their computational efficiency decreases significantly as the problem size increases, making them impractical for large-scale applications.
	
	\cite{Chen1995} formulated the multi-container, multi-box problem as a mixed-integer linear programming (MILP) model, incorporating an objective function and twelve constraints. Their approach successfully found the optimal solution for an instance with three containers and six heterogeneous boxes within 15 minutes. \cite{Martello2000} combined exact branch-and-bound algorithms with heuristic methods, achieving optimal solutions for problem instances containing up to 90 boxes within a reasonable time. \cite{Fanslau2010} proposed a tree search algorithm with region control, demonstrating promising performance at the time. More recently, \cite{Oliveira2020} formulated the 3D-BPP as an integer linear programming model, incorporating linear constraints to ensure box stability and leveraging branch-and-bound methods for optimization. \cite{Nascimento2021} combined linear programming with constraint programming to reduce the computational burden of constraint verification. Their approach effectively handled instances with up to 110 boxes of 10 different types, achieving optimal solutions for over 70\% of the cases within a time limit of 3600 seconds.
	
	\subsection{Heuristic Algorithms}
	Heuristic algorithms leverage domain-specific knowledge to establish placement rules at each step of the packing process, yielding near-optimal solutions efficiently.
	
\cite{George1980} introduced a layer-based packing strategy, where containers are filled sequentially, layer by layer. \cite{Bischoff1990} experimentally evaluated 14 heuristic rules under different box category distributions and analyzed the impact of box types on packing efficiency. \cite{Chien1998} combined layer-based methods with dynamic programming to enhance packing performance. \cite{Pisinger2002} divided the container into multiple vertical layers, selecting boxes based on a look-ahead search strategy. \cite{Eley2002} integrated greedy algorithms with tree search techniques, where the greedy approach selected boxes for optimal placement, while the tree search expanded packing sequences to explore additional solutions.
	
\cite{Bortfeldt2003} introduced a parallel tabu search algorithm for weakly heterogeneous 3D-BPP, while \cite{Eley2003} adopted a two-phase strategy, using a greedy algorithm for single-container packing followed by a heuristic approach for multiple containers. \cite{Mack2004} proposed a parallel hybrid local search algorithm combining simulated annealing and tabu search. Later advancements include the fitness-based heuristic algorithm (FDA) by \cite{He2011} and a space-blocking selection strategy by \cite{Zhu2012a}. \cite{Zhu2012} further refined these techniques by developing new block generation and selection algorithms, which were later enhanced by \cite{Araya2014, Araya2017} using beam search methods.
	
	\subsection{Metaheuristic Algorithms}
	Metaheuristic algorithms differ from heuristic approaches in incorporating randomness, generating different solutions for the same problem instance across multiple runs. This randomness allows them to explore larger solution spaces, increasing the likelihood of finding optimal or near-optimal solutions.
	
	\cite{Gehring1997} proposed a genetic algorithm (GA) that stacked boxes in tower-like formations. \cite{Bortfeldt2001} introduced a hybrid GA that divided the container into vertical layers to improve packing efficiency. \cite{Gehring2002} further developed a parallel GA specifically for highly heterogeneous 3D-BPP instances. \cite{Karabulut2004} introduced the deepest-last-left (DLL) placement strategy, integrating it with GA for enhanced performance. \cite{Moura2005} applied a greedy randomized adaptive search procedure (GRASP) to address both weakly and strongly heterogeneous 3D-BPP, considering box stability.
	
\cite{Crainic2008} proposed an extreme point-based placement rule to identify high-quality box positions efficiently. \cite{Parreno2008} introduced an overlapping space representation method to improve packing performance. \cite{Zhang2007} combined composite block generation, basic heuristic methods, and simulated annealing. \cite{Liu2011} designed a heuristic packing strategy integrated with tabu search to optimize weight distribution and stability constraints. \cite{Kang2012} improved the DLL strategy by integrating it with GA. \cite{AlvarezValdes2013} were the first to apply path relinking with GRASP for 3D-BPP, achieving promising results.
	
	Recent advancements include hybrid genetic algorithms with improved stability rules \cite{GalraoRamos2016}, simulated annealing combined with tree search \cite{Sheng2017}, and swarm optimization \cite{Zhou2017}. Saraiva et al. (2019)\cite{DiasSaraiva2019} combined metaheuristics with exact algorithms for solving single-container weakly heterogeneous 3D-BPP. \cite{Su2021} introduced a chemistry-inspired optimization metaheuristic, while \cite{Zhang2022} proposed a hybrid heuristic incorporating ant colony optimization and constructive greedy algorithms, further improving search efficiency. \cite{Huang2022} developed a differential evolution algorithm using a ternary search tree, while \cite{Fom2024} employed a non-dominated sorting genetic algorithm (NSGA) for 3D-BPP.
	
	\subsection{Deep Learning-Based Methods}
	Deep learning techniques have recently been explored to enhance the efficiency of solving 3D-BPP. \cite{Zhu2021} integrated deep learning with tree search algorithms to optimize packing strategies. \cite{Que2023} employed a combination of Transformer models and reinforcement learning to address 3D-BPP, while \cite{Chien2024} leveraged deep learning to guide genetic algorithms, improving local search effectiveness.
	
	The 3D bin packing problem remains a challenging combinatorial optimization problem. Exact algorithms provide optimal solutions but are computationally prohibitive for large-scale instances. Heuristic and metaheuristic approaches have been widely applied due to their efficiency and adaptability. More recently, deep learning-based methods have demonstrated promise in improving solution quality and computational efficiency. Since heuristic methods can consider both efficiency and stability, the method employed in this study is heuristic.

	\section{Methodology}
	In this study, we organize the nuclear waste boxes into large blocks and employ a beam search algorithm to determine an optimal loading scheme. Using blocks helps reduce the solution space, thus enhancing the efficiency of the search process. As demonstrated by several previous studies\cite{Araya2014, Araya2017}, loading using blocks can yield promising results. The beam search algorithm, widely recognized for its effectiveness in combinatorial optimization problems, has proven successful in three-dimensional and two-dimensional bin-packing problems. For these reasons, we chose to apply the beam search algorithm in this work.
	
	Prior research \cite{Zhu2012} on three-dimensional block bin packing problem identified six key elements for an effective algorithm: 
\begin{enumerate}[(1)]
	\item Representation of the remaining space in the container.  
	\item Block generation, i.e., the block generation algorithm.  
	\item Selection of remaining space, i.e., the remaining space selection algorithm.  
	\item Selection of blocks, i.e., the block selection algorithm.  
	\item Placement of the selected blocks into the chosen space and updating the space accordingly.  
	\item The search strategy, i.e., the search algorithm.
\end{enumerate}
	
	In this paper, we will introduce our approach based on these six aspects and comprehensively explain how each element is addressed in our method.
	\subsection{Spatial Representation}
	In three-dimensional block packing, each remaining space is represented as a rectangular prism. When a block is placed into a given space, it divides the space into three smaller subspaces. Depending on whether these newly created subspaces overlap, the representation of the remaining space can be classified into two types: partial space representation and overlapping space representation.

\begin{figure}
		\centering
		\includegraphics[width=.6\columnwidth]{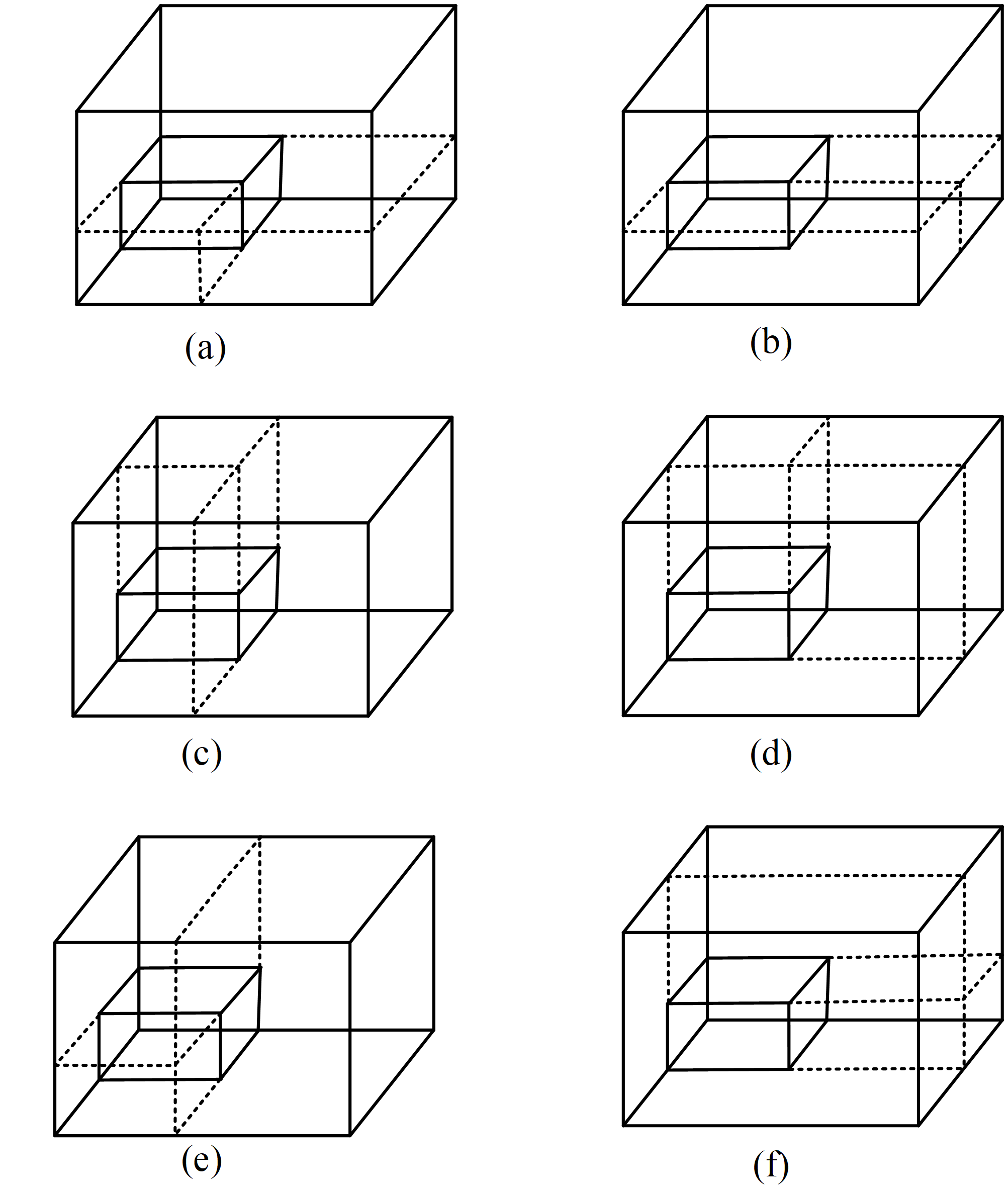}
		\caption{partial space}
		\label{paspace}
	\end{figure}

The partial space representation divides the remaining space into three non-overlapping rectangular prisms. Based on different division methods, the remaining space has six possible configurations, as illustrated in Figure \ref{paspace}. The majority of studies employing partial representation typically use the division methods shown in (c) or (d) of Figure \ref{paspace}. Some studies, such as the one by \cite{Zhang2007}, combine these two division methods.
	
\begin{figure}
		\centering
		\includegraphics[width=.9\columnwidth]{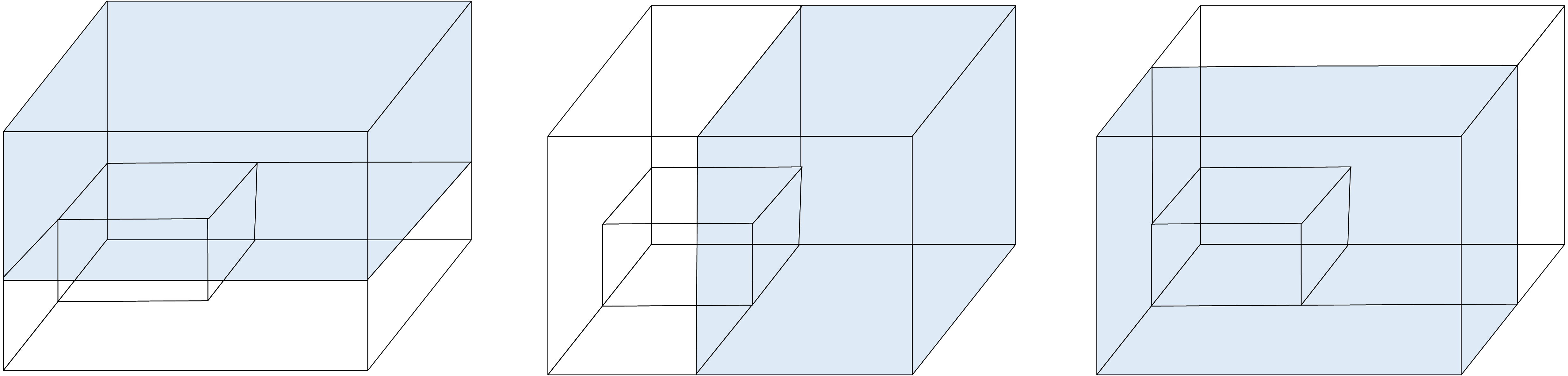}
		\caption{overlapping space.}
		\label{cospace}
	\end{figure}

	On the other hand, the overlapping space representation depicts the remaining space as three overlapping rectangular spaces, as shown in Figure \ref{cospace}. For clarity, three separate diagrams are provided in the figure to display the three subspaces resulting from the overlapping representation.
	
	In this study, the overlapping space representation is preferred, as it can capture more detailed information about the remaining space, which typically leads to more efficient packing outcomes.
	
	\subsection{Block Generation Algorithm}
	The block generation algorithm combines the boxes to be packed into larger blocks according to specific rules, while individual boxes can also form independent blocks.
In this study, we first combine nuclear waste boxes of the same size into blocks, forming simple blocks, where individual nuclear waste boxes can also form independent blocks. Then, simple blocks are paired and combined to create larger blocks called complex blocks. All combined blocks must satisfy the following conditions: 
	\begin{enumerate}
		\item The number of nuclear waste boxes must be sufficient to form the block.
		\item The block's length must be less than or equal to the length of the disposal pool, and its width must be less than the width of the disposal pool.
	\end{enumerate}

	To determine whether two blocks, \( b_1 \) and \( b_2 \), can be combined into a larger block, we follow the rules established in previous research:
	\begin{enumerate}
		\item When combining along the x-direction, if \( b_1.\text{width} \times b_1.\text{height} \geq b_2.\text{width} \times b_2.\text{height} \), block \( b_2 \) can only be placed to the right of \( b_1 \).
		\item When combining along the y-direction, if \( b_1.\text{length} \times b_1.\text{height} \geq b_2.\text{length} \times b_2.\text{height} \), block \( b_2 \) can only be placed in front of \( b_1 \).
		\item When combining along the z-direction, if \( b_1.\text{length} \times b_1.\text{width} \geq b_2.\text{length} \times b_2.\text{width} \), block \( b_2 \) can only be placed above \( b_1 \).
	\end{enumerate}

	Let \( b_1 \) and \( b_2 \) combined form block \( b_3 \). Additionally, block \( b_3 \) must satisfy the condition that the volume of the boxes inside \( b_3 \), divided by the total volume of \( b_3 \), is greater than or equal to the minimum volume utilization rate \( \text{minFillRate} \).
	
\begin{figure}
		\centering
		\includegraphics[width=.9\columnwidth]{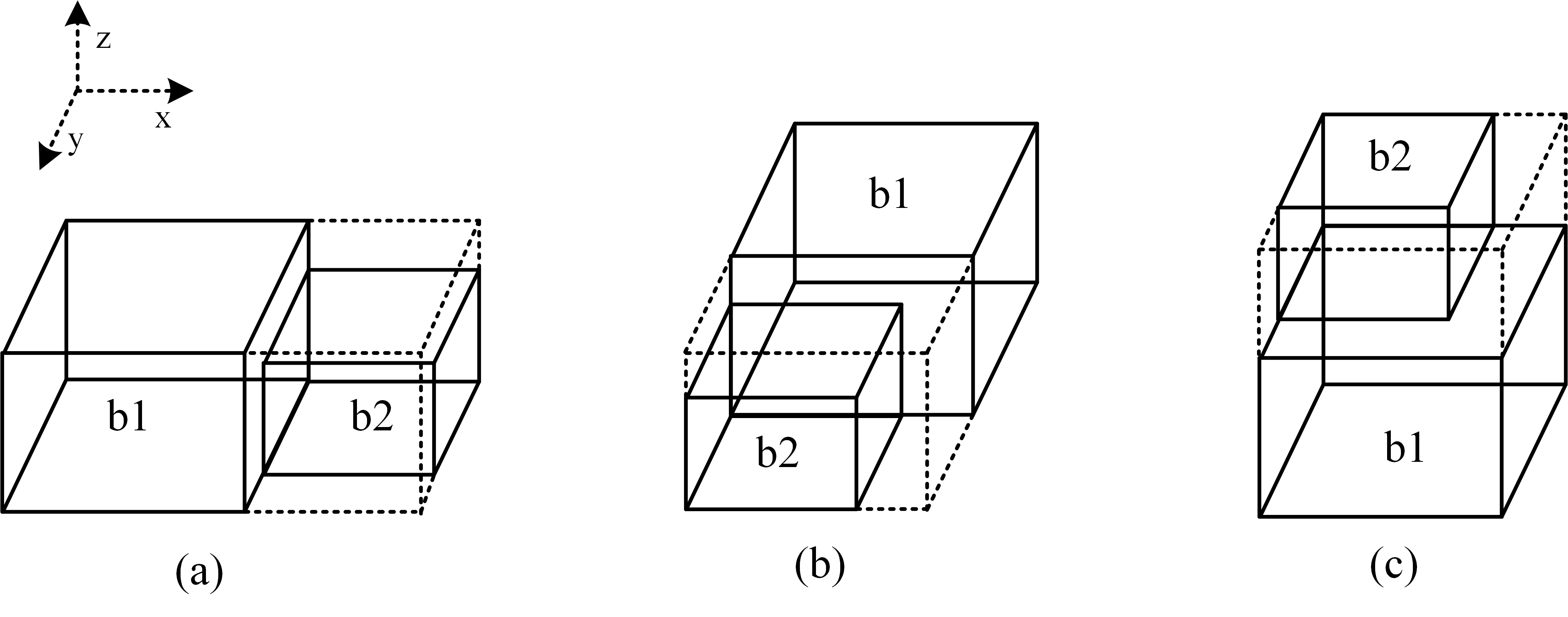}
		\caption{Block connection in three directions.}
		\label{cb}
	\end{figure}

	Figure \ref{cb} illustrates the results of combining blocks \( b_1 \) and \( b_2 \) along the three directions.
	The process of generating blocks used in this study is shown in Algorithm\ref{blockg}. Furthermore, the dose rate at a distance of one meter from the center of the block is equivalent to the sum of the dose rates at a distance of one meter from the centers of the individual boxes contained within the block.

	\begin{algorithm}[ht]
		\caption{Block Generation}
		\label{blockg}
		\SetAlgoLined
		\KwIn{$boxList$, $maxNum$, $minFillRate$}
		\KwOut{$blockList$}
		
		$blockList \gets$ generate simple blocks\;  
		$pList \gets blockList$\;
		
		\While{ $|blockList| < maxNum$}{
			$newBlockList \gets \emptyset$\;
			\ForEach{$block_1 \in pList$}{
				\ForEach{$block_2 \in blockList$}{
					\ForEach{direction $\in$ \{x, y, z\}}{
						Concatenate $block_1$ and $block_2$ based on the conditions, and add the generated block to $newBlockList$.
					}
				}
			}
			
			\ForEach{$block \in newBlockList$}{
				\If{$block \notin blockList$ \textbf{and} $|blockList| < maxNum$}{
					$blockList \gets blockList \cup \{block\}$\;
				}
			}
			
			$pList \gets newBlockList$\;
		}
		\Return{$blockList$}\;
	\end{algorithm}

	\subsection{Spatial Selection Algorithm}
	One of the main objectives of this study is to minimize the total dose rate of all the waste boxes in the disposal pool. To achieve this, we select the space with the lowest position and the farthest left within the remaining available spaces at each step. This approach is based on the fact that the dose rate of each nuclear waste box is inversely proportional to the vertical distance between the center of the box and the top of the disposal pool. By selecting the space with the lowest position, we ensure that boxes with higher activity values are placed as deep as possible within the disposal pool. This placement strategy effectively reduces the overall sum of the dose rates for all the nuclear waste boxes in the pool, as it minimizes the exposure of high-activity boxes to the pool's surface.
	\subsection{Block Selection Algorithm}
	For a given space \( s \) and any block \( b_i \) from the \( \text{blockList} \), we first identify all the blocks that can be placed within space \( s \). Then, the scores of these blocks are computed individually. Finally, the top-ranked blocks are selected based on their scores and ordered in descending order. The formula for calculating the block score is as (\ref{fbs}).
	\begin{equation}\label{fbs}
		f(b_i, s) = \frac{V(b_i) - V_{\text{loss}}(b_i, s) - V_{\text{waste}}(b_i)}{V(s)} +\alpha \times A_{nor}(b_i)
	\end{equation}
	In this formula, \( V(b_i) \) represents the volume of block \( b_i \), \( V_{\text{loss}}(b_i, s) \) denotes the wasted volume in space \( s \) after placing block \( b_i \), and \( V_{\text{waste}}(b_i) \) is the volume of block \( b_i \) minus the volume occupied by the boxes within it. \( V(s) \) refers to the total volume of space \( s \), \( A_{\text{nor}}(b_i) \) is the min-max normalized value of the dose rate at a distance of one meter from the center of the block \( b_i \), and $\alpha$ is a parameter.
	
	For \( V_{\text{loss}}(b_i, s) \), we define \( l(b_i) \), \( w(b_i) \), and \( h(b_i) \) as the length, width, and height of block \( b_i \), respectively, and \( l(s) \), \( w(s) \), and \( h(s) \) as the length, width, and height of space \( s \), respectively. The maximum linear combination of the lengths, widths, and heights of the remaining boxes that do not exceed the corresponding dimensions of the available space are represented by \( l_{\text{max}} \), \( w_{\text{max}} \), and \( h_{\text{max}} \). Specifically, \( l_{\text{max}} \) is the maximum linear combination of the lengths of all remaining boxes that does not exceed \( l(s) - l(b_i) \), \( w_{\text{max}} \) is the maximum linear combination of the widths of all remaining boxes that does not exceed \( w(s) - w(b_i) \), and \( h_{\text{max}} \) is the maximum linear combination of the heights of all remaining boxes that does not exceed \( h(s) - h(b_i) \). The final formula for calculating \( V_{\text{loss}}(b_i, s) \) is as (\ref{vloss}).
	\begin{equation}\label{vloss}
		\begin{split}
			V_{\text{loss}}(b_i, s) = & V(s) - (l(b_i) + l_{\text{max}}) \cdot  \\
			& (w(b_i) + w_{\text{max}}) \cdot (h(b_i) + h_{\text{max}})
		\end{split}
	\end{equation}
	
	\subsection{Space Update}
\begin{figure}
		\centering
		\includegraphics[width=.7\columnwidth]{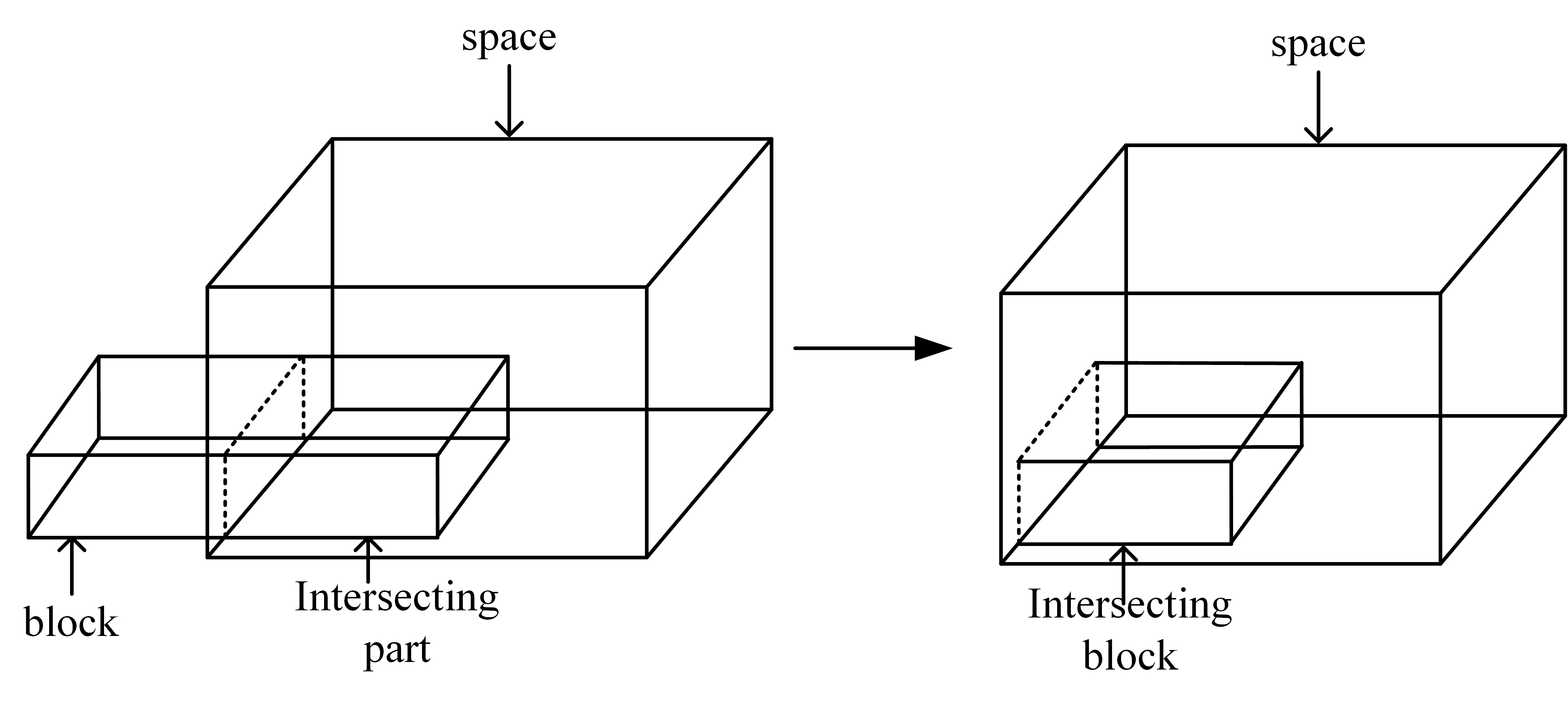}
		\caption{Space update.}
		\label{renewspace}
	\end{figure}

	The selected block is positioned at the bottom-left corner of the chosen space. The space update procedure is as follows: first, the selected space is removed from the remaining spaces. Then, the new spaces created as a result of placing the block are calculated. Finally, all spaces that overlap with the block are updated. The update method for these overlapping spaces involves identifying the intersecting region between the block and each space. This overlapping region is treated as a new block, and the space is updated accordingly. A schematic representation of this update process is provided in Figure \ref{renewspace}.
	
	\subsection{Search Algorithm}
	We employ a beam search algorithm to explore potential placement schemes. The framework of the algorithm is outlined in Algorithm \ref{beam}.

Before initiating the search, we define the algorithm's runtime, \( t \), and generate the initial state, \( stateInit \). This state includes the number of each type of nuclear waste box remaining \( remainbox \), a list of blocks \( blockList \), a list of available spaces \( spaceList \), and a state score. Initially, \( remainbox \) represents the total quantity of each type of nuclear waste box yet to be placed, \( blockList \) is the list of blocks generated by the block generation algorithm, and \( spaceList \) contains the entire disposal pool.
	\begin{algorithm}[ht]
		\caption{Beam search}
		\label{beam}
		\KwIn{$stateInit$, $bestState_a$}
		\KwOut{$bestState_a$}
		
		$w \gets 1$ \;
		
		\While{Time Limit Is Reached}{
			$stateList \gets \emptyset$ \;
			Add $stateInit$ to $stateList$ \;
			
			\While{$stateList$ is not empty}{
				$successors \gets \emptyset$ \;
				\ForEach{$state \in stateList$}{
					\If{$state$ equals $stateInit$}{
						$succ \gets$ expand($state$, $w \times w$) \;
					}
					\Else{
						$succ \gets$ expand($state$, $w$) \;
					}
					Add all states in $succ$ to $successors$ \;
				}
				
				$stateList \gets \emptyset$ \;
				
				\If{$successors$ is not empty}{
					\ForEach{$state \in successors$}{
						greedy($state$, $bestState_a$) \;
					}
				}
				Add the top $w$ states in $successors$ ranked by score to $stateList$ \;
			}
			
			$w \gets \sqrt{2} \times w$ \;
		}
		
		\Return{$bestState_a$} \;
	\end{algorithm}
	
	During the search process, for each state, the space selection algorithm is first applied to choose a space for block placement. Then, the block selection method is used to select \( w \) blocks from \( blockList \), which are placed in the chosen space, generating \( w \) new states; this process is shown as Algorithm \ref{expand}. For each of these new states, a greedy method is employed to calculate the corresponding score, and the best solution encountered during the greedy is recorded. Subsequently, the higher-score states are selected as nodes for the next expansion. After each search round, the width of the search, \( w \), is updated by multiplying it by \( \sqrt{2} \).
	\begin{algorithm}[ht]
		\caption{expand}
		\label{expand}
		\KwIn{$state$, $w$}
		\KwOut{$succ$}
		
		$blocks \gets \emptyset$ \;
		
		\While{$blocks = \emptyset$ \textbf{and} $spaceList \neq \emptyset$}{
			$space \gets$ select a space from $spaceList$\;
			$blocks \gets$ select w block from $blockList$\;
			
			\If{$blocks$ is empty}{
				Remove $space$ from $spaceList$ \;
			}
			\Else{
				\textbf{break}
			}
		}
		
		$succ \gets \emptyset$ \;
		
		\ForEach{$block \in blocks$}{
			Place the $block$ into $space$ to generate a new state, and then add the new state into $succ$ \;
		}
		
		\Return{$succ$} \;
	\end{algorithm}
	
For the generation of new states, the algorithm first removes blocks that can no longer be formed with the remaining nuclear waste boxes from \( blockList \). Then, the space update method is applied to update \( spaceList \), and the number of remaining nuclear waste boxes, \( remainbox \), is adjusted accordingly.

The greedy algorithm, as shown in algorithm \ref{greedy}, demonstrates the process for calculating each state's score. A space is selected for each state, and a block is chosen for placement in that space. This process is repeated until no space remains. Once all possible placements are exhausted, the current placement scheme is compared with the optimal solution. The scheme with the higher space utilization and lower dose rate is updated as the optimal solution. Additionally, the loading volume is used as the score for each state.
	
	\begin{algorithm}[ht]
		\SetAlgoLined
		\KwIn{$state$, $bestState$}
		
		\While{state.spaceList is not empty}{
			$space \gets$ select a space from $spaceList$\;
			$block \gets$ select a block from $blockList$\;
			
			\eIf{$block$ is not empty}{
				Place $block$ into $space$ and renew $state$\;
			}{
				remove $space$ from $spaceList$ \;
			}
		}
		
		Compare $state$ and $bestState$, and store the one with higher space utilization and lower dose rate in $bestState$ \;
		
		Set the score of $state$ to the loading volume \;
		
		\caption{greedy}
		\label{greedy}
	\end{algorithm}
	
	\section{Experiment}
	In this section, we describe the datasets employed in this study. Using these datasets, we evaluate the performance of the proposed algorithm by testing the utilization rate and dose rate under different values of parameter $\alpha$.
	
	\subsection{Datesets}
	We constructed 16 datasets, labeled set1 to set16, each containing 100 nuclear waste packing problems. These datasets were developed using real-world nuclear waste packing data from nuclear power plants, in combination with the three-dimensional bin packing problem datasets provided by \cite{Bischoff1995} and \cite{Davies1999}. A summary of the number of box types, disposal pool sizes, and the number of boxes per problem for each dataset is presented in Table \ref{datainfo}.
	\begin{table}[width=.9\linewidth,cols=5,pos=h]
		\caption{Datasets info.}\label{datainfo}
		\begin{tabular*}{\tblwidth}{@{} CCCCC@{} }
			\toprule
			\makecell{Datasets} & \makecell{Box type\\ count} & \makecell{Pool size\\ (L,W,H)} & \makecell{No. of \\ problems} & \makecell{Boxes per \\ problem}  \\ 
			set1 & 1 & 587 233 220 & 100 & 205  \\ 
			set2 & 3 & 587 233 220 & 100 & 150  \\ 
			set3 & 5 & 587 233 220 & 100 & 136  \\ 
			set4 & 8 & 587 233 220 & 100 & 134  \\ 
			set5 & 10 & 587 233 220 & 100 & 132  \\ 
			set6 & 12 & 587 233 220 & 100 & 132  \\ 
			set7 & 15 & 587 233 220 & 100 & 131  \\ 
			set8 & 20 & 587 233 220 & 100 & 130  \\ 
			set9 & 30 & 587 233 220 & 100 & 130  \\ 
			set10 & 40 & 587 233 220 & 100 & 128  \\ 
			set11 & 50 & 587 233 220 & 100 & 130  \\ 
			set12 & 60 & 587 233 220 & 100 & 129  \\ 
			set13 & 70 & 587 233 220 & 100 & 130  \\ 
			set14 & 80 & 587 233 220 & 100 & 130  \\ 
			set15 & 90 & 587 233 220 & 100 & 129  \\ 
			set16 & 100 & 587 233 220 & 100 & 129  \\ 
			\bottomrule
		\end{tabular*}
	\end{table}
	
	\subsection{Experimental Results}
	
\begin{figure}
		\centering
		\includegraphics[width=.9\columnwidth]{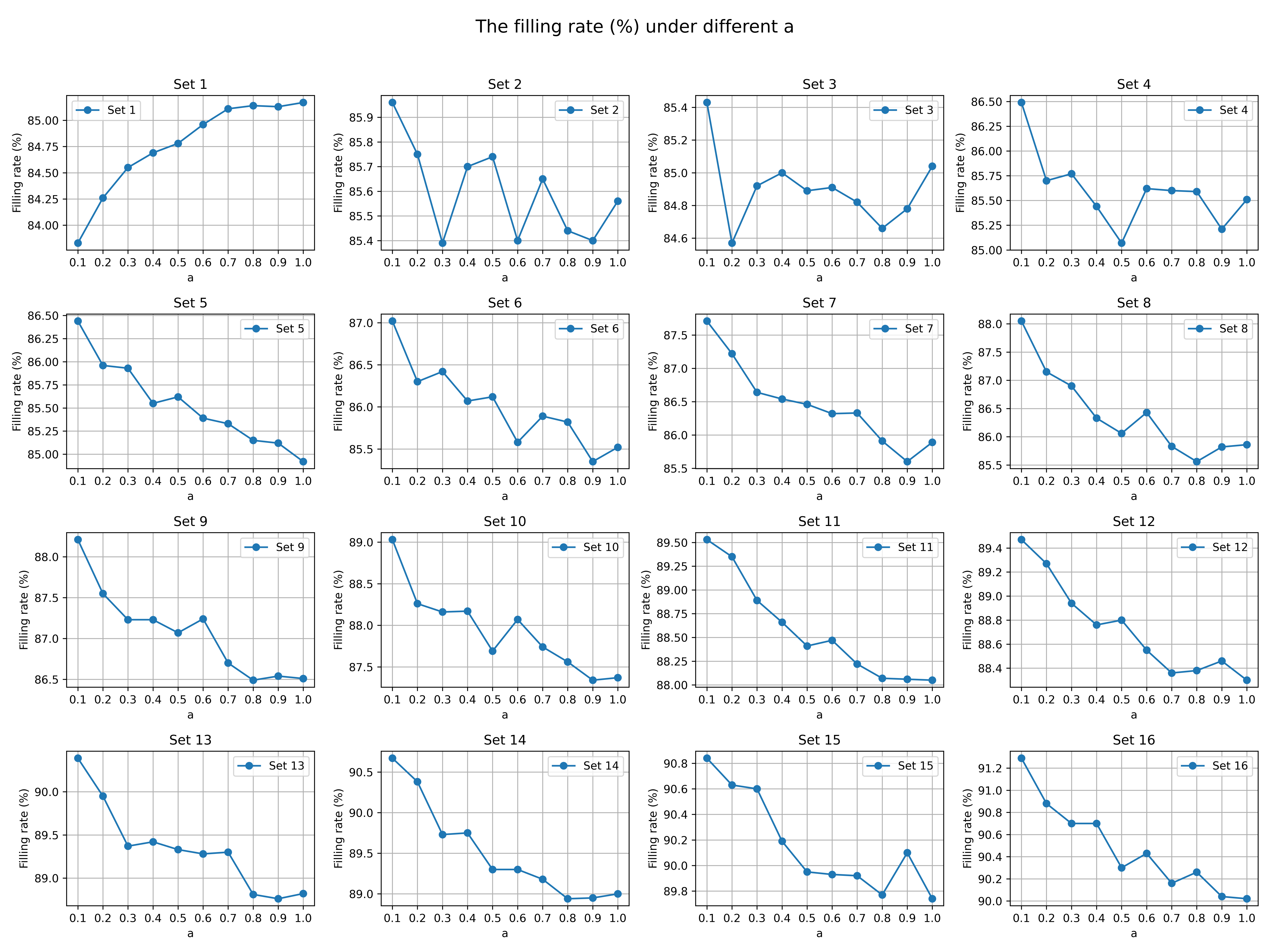}
		\caption{Spatial utilization rate (\%) under varying $\alpha$.}
		\label{fillratepng}
	\end{figure}

	Table \ref{fillrate} and Figure \ref{fillratepng} present the space utilization of the BSNA for different values of parameter \(\alpha\), with a fixed running time of \(t = 30s\). The filling rate decreases for datasets set1 through set16 as \(\alpha\) increases. This occurs because a higher value of a raises the weight of \(A_{nor}\) in the block scoring formula, shifting the selection criteria from block volume to dose rate and thus reducing space utilization. In contrast, for dataset set1—where all boxes are of identical size—the volume-based term remains constant, effectively reducing the scoring formula to one that emphasizes dose rate as shown in \ref{fbs2}, which results in an increased filling rate as \(\alpha\) increases.
	
	\begin{equation}\label{fbs2}
		f(b_i, s) = \alpha \times A_{nor}(b_i)
	\end{equation}
	
\begin{figure}
		\centering
		\includegraphics[width=.9\columnwidth]{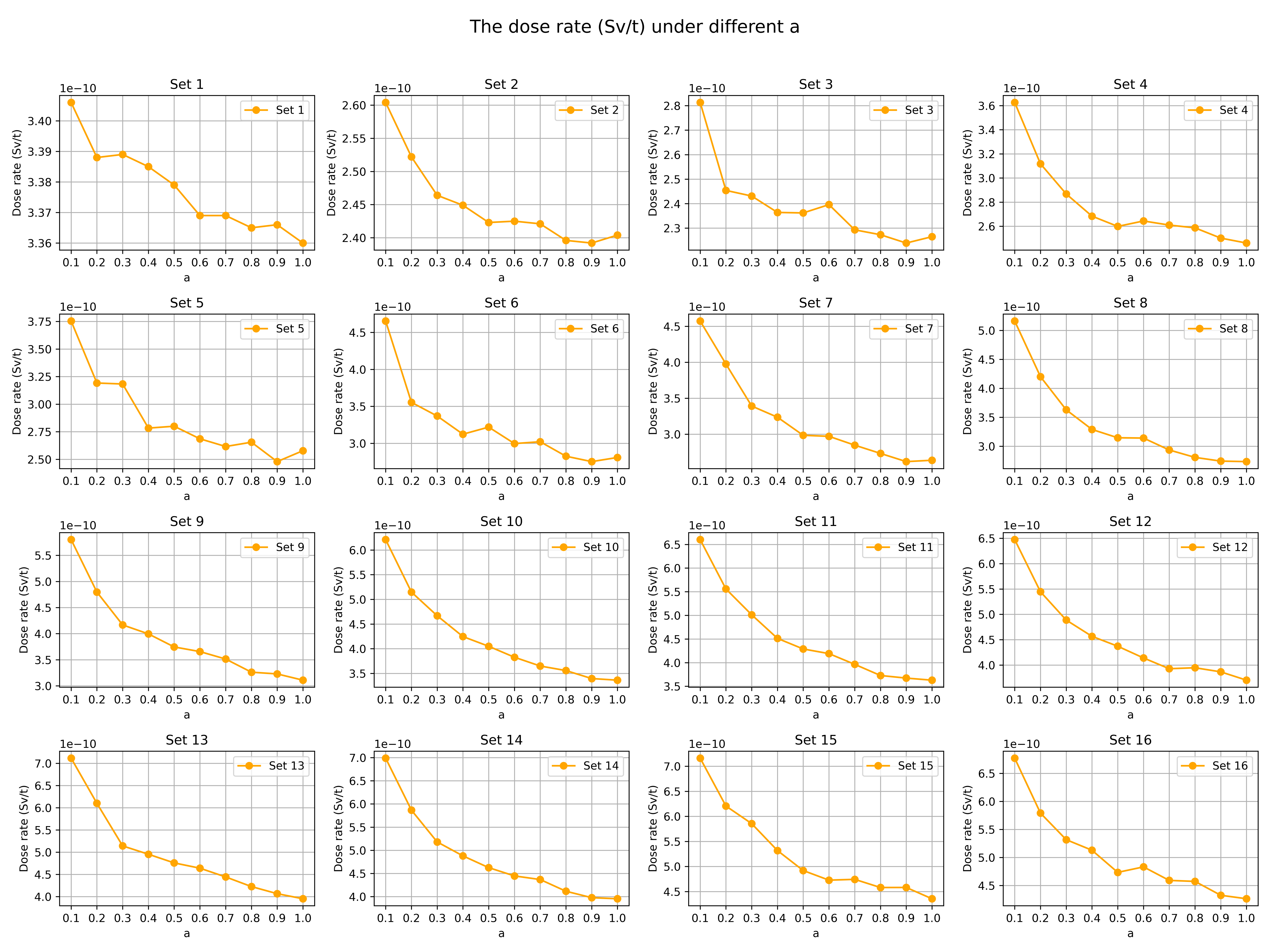}
		\caption{Dose rate (Sv/t) under different \(\alpha\).}
		\label{doseratepng}
	\end{figure}

	Table \ref{doserate} and Figure \ref{doseratepng} present the dose rate results of the BSNA for different values of parameter \(\alpha\), with a fixed running time of  \(t = 30s\). The findings indicate that increasing \(\alpha\) reduces the overall dose rate across all datasets. This trend arises because a higher \(\alpha\) increases the weight of \(A_{nor}\) in the block scoring formula. Consequently, the algorithm preferentially selects blocks with higher dose rates and assigns them to the lowest available positions, thereby minimizing the cumulative dose rate in the disposal pool.
	
\begin{figure}
		\centering
		\includegraphics[width=.9\columnwidth]{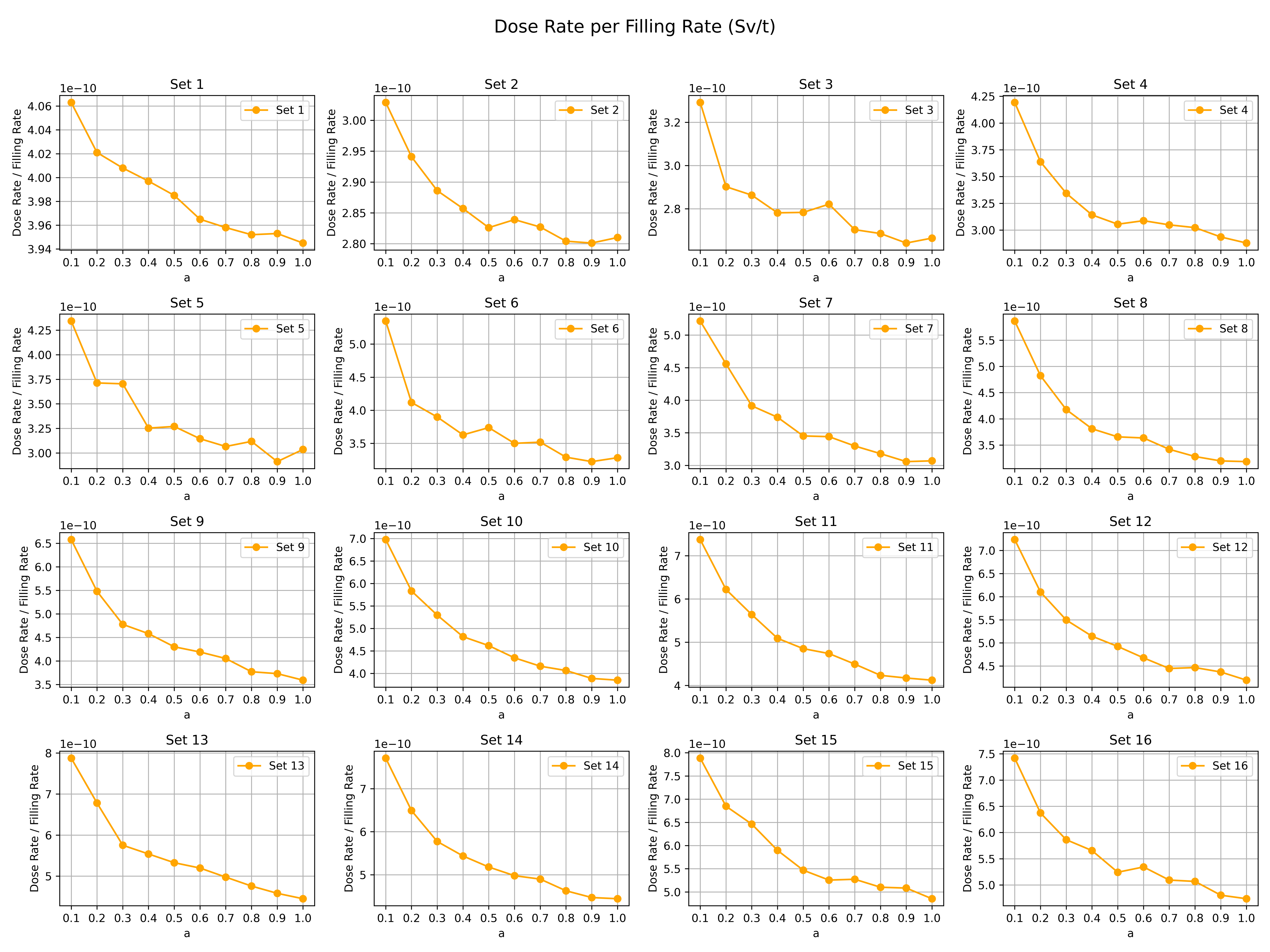}
		\caption{Dose rate per filling rate (Sv/t).}
		\label{dperfpng}
	\end{figure}

One of the primary objectives of the proposed algorithm is to maximize the disposal pool's space utilization. However, the results in Table \ref{fillrate} seem to contradict this aim. To clarify the relationship between space utilization and dose rate, Table \ref{dperftable} and Figure \ref{dperfpng} present the ratio of dose rate to space utilization rate—an indicator of the dose rate per unit volume of the waste boxes. The results clearly show that as parameter \(\alpha\) increases, the dose rate per unit volume decreases, and indicate that increasing parameter \(\alpha\) leads the algorithm to select boxes with a lower dose rate per unit volume, potentially allowing for a higher loading rate without an increase in the overall dose rate. To further validate this observation, we set \(\alpha\) to 0.6 and evaluated the spatial utilization and dose rate across all datasets at running times of 30 s, 60 s, 90 s, and 120 s, as shown in Table \ref{fortime}. The results demonstrate that extending the search time with a fixed parameter \(\alpha = 0.6\) improves spatial utilization while concurrently reducing the dose rate.
	
	In summary, for a fixed search time, increasing parameter \(\alpha\) reduces both the dose and filling rates. Conversely, for a fixed value of \(\alpha\), extending the algorithm’s search time decreases the dose rate while enhancing the filling rate. This indicates that, within a limited timeframe, increasing \(\alpha\) effectively reduces the disposal pool’s dose rate, whereas longer search times improve space utilization.
	
	\begin{table*}[width=2\linewidth,cols=11,pos=h]
		\caption{The spatial utilization rate (\%) under varying $\alpha$ values.(\(t = 30s\))}\label{fillrate}
		\begin{tabular*}{\tblwidth}{@{} LCCCCCCCCCC@{} }
			\toprule
			\multirow{2}{*}{Dataset} & \multicolumn{10}{c}{$\alpha$} \\
			\cmidrule{2-11}
			& 0.1 & 0.2 & 0.3 & 0.4 & 0.5 & 0.69 & 0.7 & 0.8 & 0.9 &1.0 \\
			\midrule
			set1 & 83.83 & 84.26 & 84.55 & 84.69 & 84.78 & 84.96 & 85.11 & 85.14 & 85.13 & 85.17  \\ 
			set2 & 85.96 & 85.75 & 85.39 & 85.70 & 85.74 & 85.40 & 85.65 & 85.44 & 85.40 & 85.56  \\ 
			set3 & 85.43 & 84.57 & 84.92 & 85.00 & 84.89 & 84.91 & 84.82 & 84.66 & 84.78 & 85.04  \\ 
			set4 & 86.49 & 85.70 & 85.77 & 85.44 & 85.07 & 85.62 & 85.60 & 85.59 & 85.21 & 85.51  \\ 
			set5 & 86.44 & 85.96 & 85.93 & 85.55 & 85.62 & 85.39 & 85.33 & 85.15 & 85.12 & 84.92  \\ 
			set6 & 87.02 & 86.30 & 86.42 & 86.07 & 86.12 & 85.58 & 85.89 & 85.82 & 85.35 & 85.52  \\ 
			set7 & 87.71 & 87.22 & 86.64 & 86.54 & 86.46 & 86.32 & 86.33 & 85.91 & 85.60 & 85.89  \\ 
			set8 & 88.05 & 87.15 & 86.90 & 86.33 & 86.06 & 86.43 & 85.83 & 85.56 & 85.82 & 85.86  \\ 
			set9 & 88.21 & 87.55 & 87.23 & 87.23 & 87.07 & 87.24 & 86.70 & 86.49 & 86.54 & 86.51  \\ 
			set10 & 89.03 & 88.26 & 88.16 & 88.17 & 87.69 & 88.07 & 87.74 & 87.56 & 87.34 & 87.37  \\ 
			set11 & 89.53 & 89.35 & 88.89 & 88.66 & 88.41 & 88.47 & 88.22 & 88.07 & 88.06 & 88.05  \\ 
			set12 & 89.47 & 89.27 & 88.94 & 88.76 & 88.80 & 88.55 & 88.36 & 88.38 & 88.46 & 88.30  \\ 
			set13 & 90.39 & 89.95 & 89.37 & 89.42 & 89.33 & 89.28 & 89.30 & 88.81 & 88.76 & 88.82  \\ 
			set14 & 90.67 & 90.38 & 89.73 & 89.75 & 89.30 & 89.30 & 89.18 & 88.94 & 88.95 & 89.00  \\ 
			set15 & 90.84 & 90.63 & 90.60 & 90.19 & 89.95 & 89.93 & 89.92 & 89.77 & 90.10 & 89.74  \\ 
			set16 & 91.29 & 90.88 & 90.70 & 90.70 & 90.30 & 90.43 & 90.16 & 90.26 & 90.04 & 90.02  \\ 
			
			\bottomrule
		\end{tabular*}
	\end{table*}
	
	\begin{table*}[width=2\linewidth,cols=11,pos=h]
		\caption{The dose rate (Sv/t) under different \(\alpha\) values.(\(t = 30s\))}\label{doserate}
		\begin{tabular*}{\tblwidth}{@{} LCCCCCCCCCC@{} }
			\toprule
			\multirow{2}{*}{Dataset} & \multicolumn{10}{c}{$\alpha$} \\
			\cmidrule{2-11}
			& 0.1 & 0.2 & 0.3 & 0.4 & 0.5 & 0.6 & 0.7 & 0.8 & 0.9 &1.0 \\
			\midrule
			set1 & 3.406E-10 & 3.388E-10 & 3.389E-10 & 3.385E-10 & 3.379E-10 & 3.369E-10 & 3.369E-10 & 3.365E-10 & 3.366E-10 & 3.360E-10  \\ 
			set2 & 2.604E-10 & 2.522E-10 & 2.464E-10 & 2.449E-10 & 2.423E-10 & 2.425E-10 & 2.421E-10 & 2.396E-10 & 2.392E-10 & 2.404E-10  \\ 
			set3 & 2.813E-10 & 2.454E-10 & 2.431E-10 & 2.364E-10 & 2.362E-10 & 2.396E-10 & 2.293E-10 & 2.273E-10 & 2.239E-10 & 2.265E-10  \\ 
			set4 & 3.626E-10 & 3.118E-10 & 2.868E-10 & 2.684E-10 & 2.599E-10 & 2.644E-10 & 2.610E-10 & 2.587E-10 & 2.502E-10 & 2.461E-10  \\ 
			set5 & 3.754E-10 & 3.192E-10 & 3.183E-10 & 2.783E-10 & 2.800E-10 & 2.686E-10 & 2.616E-10 & 2.655E-10 & 2.479E-10 & 2.578E-10  \\ 
			set6 & 4.654E-10 & 3.554E-10 & 3.369E-10 & 3.123E-10 & 3.219E-10 & 2.996E-10 & 3.021E-10 & 2.825E-10 & 2.752E-10 & 2.808E-10  \\ 
			set7 & 4.574E-10 & 3.976E-10 & 3.391E-10 & 3.236E-10 & 2.984E-10 & 2.970E-10 & 2.847E-10 & 2.732E-10 & 2.618E-10 & 2.637E-10  \\ 
			set8 & 5.163E-10 & 4.202E-10 & 3.629E-10 & 3.290E-10 & 3.146E-10 & 3.141E-10 & 2.935E-10 & 2.807E-10 & 2.743E-10 & 2.734E-10  \\ 
			set9 & 5.802E-10 & 4.797E-10 & 4.167E-10 & 3.995E-10 & 3.745E-10 & 3.656E-10 & 3.513E-10 & 3.261E-10 & 3.228E-10 & 3.108E-10  \\ 
			set10 & 6.211E-10 & 5.148E-10 & 4.668E-10 & 4.246E-10 & 4.049E-10 & 3.827E-10 & 3.649E-10 & 3.557E-10 & 3.395E-10 & 3.362E-10  \\ 
			set11 & 6.602E-10 & 5.558E-10 & 5.012E-10 & 4.513E-10 & 4.289E-10 & 4.191E-10 & 3.964E-10 & 3.728E-10 & 3.673E-10 & 3.629E-10  \\ 
			set12 & 6.473E-10 & 5.446E-10 & 4.888E-10 & 4.567E-10 & 4.372E-10 & 4.141E-10 & 3.929E-10 & 3.948E-10 & 3.867E-10 & 3.703E-10  \\ 
			set13 & 7.117E-10 & 6.101E-10 & 5.142E-10 & 4.954E-10 & 4.759E-10 & 4.640E-10 & 4.443E-10 & 4.226E-10 & 4.066E-10 & 3.952E-10  \\ 
			set14 & 6.989E-10 & 5.867E-10 & 5.180E-10 & 4.879E-10 & 4.625E-10 & 4.446E-10 & 4.368E-10 & 4.115E-10 & 3.976E-10 & 3.953E-10  \\ 
			set15 & 7.161E-10 & 6.206E-10 & 5.855E-10 & 5.318E-10 & 4.919E-10 & 4.727E-10 & 4.743E-10 & 4.580E-10 & 4.581E-10 & 4.356E-10  \\ 
			set16 & 6.771E-10 & 5.792E-10 & 5.318E-10 & 5.130E-10 & 4.735E-10 & 4.833E-10 & 4.592E-10 & 4.574E-10 & 4.327E-10 & 4.265E-10  \\ 
			
			\bottomrule
		\end{tabular*}
	\end{table*}
	
	\begin{table*}[width=2\linewidth,cols=11,pos=h]
		\caption{The dose rate of unit space utilization under different \(\alpha\) values.}\label{dperftable}
		\begin{tabular*}{\tblwidth}{@{} LCCCCCCCCCC@{} }
			\toprule
			\multirow{2}{*}{Dataset} & \multicolumn{10}{c}{$\alpha$} \\
			\cmidrule{2-11}
			& 0.1 & 0.2 & 0.3 & 0.4 & 0.5 & 0.6 & 0.7 & 0.8 & 0.9 &1.0 \\
			\midrule
			set1 & 4.063E-10 & 4.021E-10 & 4.008E-10 & 3.997E-10 & 3.985E-10 & 3.965E-10 & 3.958E-10 & 3.952E-10 & 3.953E-10 & 3.945E-10  \\ 
			set2 & 3.029E-10 & 2.941E-10 & 2.886E-10 & 2.857E-10 & 2.826E-10 & 2.839E-10 & 2.827E-10 & 2.804E-10 & 2.801E-10 & 2.810E-10  \\ 
			set3 & 3.292E-10 & 2.902E-10 & 2.863E-10 & 2.781E-10 & 2.783E-10 & 2.821E-10 & 2.703E-10 & 2.685E-10 & 2.641E-10 & 2.664E-10  \\ 
			set4 & 4.192E-10 & 3.638E-10 & 3.344E-10 & 3.142E-10 & 3.055E-10 & 3.088E-10 & 3.049E-10 & 3.023E-10 & 2.936E-10 & 2.879E-10  \\ 
			set5 & 4.343E-10 & 3.713E-10 & 3.704E-10 & 3.253E-10 & 3.270E-10 & 3.145E-10 & 3.066E-10 & 3.118E-10 & 2.912E-10 & 3.036E-10  \\ 
			set6 & 5.348E-10 & 4.118E-10 & 3.898E-10 & 3.628E-10 & 3.738E-10 & 3.500E-10 & 3.518E-10 & 3.291E-10 & 3.224E-10 & 3.283E-10  \\ 
			set7 & 5.215E-10 & 4.558E-10 & 3.914E-10 & 3.739E-10 & 3.452E-10 & 3.441E-10 & 3.298E-10 & 3.180E-10 & 3.058E-10 & 3.070E-10  \\ 
			set8 & 5.864E-10 & 4.822E-10 & 4.176E-10 & 3.810E-10 & 3.656E-10 & 3.635E-10 & 3.419E-10 & 3.281E-10 & 3.197E-10 & 3.184E-10  \\ 
			set9 & 6.577E-10 & 5.480E-10 & 4.777E-10 & 4.580E-10 & 4.301E-10 & 4.190E-10 & 4.052E-10 & 3.771E-10 & 3.731E-10 & 3.592E-10  \\ 
			set10 & 6.976E-10 & 5.833E-10 & 5.295E-10 & 4.816E-10 & 4.617E-10 & 4.345E-10 & 4.159E-10 & 4.063E-10 & 3.887E-10 & 3.848E-10  \\ 
			set11 & 7.374E-10 & 6.221E-10 & 5.639E-10 & 5.090E-10 & 4.851E-10 & 4.737E-10 & 4.493E-10 & 4.233E-10 & 4.171E-10 & 4.121E-10  \\ 
			set12 & 7.235E-10 & 6.101E-10 & 5.496E-10 & 5.145E-10 & 4.924E-10 & 4.677E-10 & 4.447E-10 & 4.467E-10 & 4.371E-10 & 4.194E-10  \\ 
			set13 & 7.874E-10 & 6.783E-10 & 5.753E-10 & 5.540E-10 & 5.327E-10 & 5.197E-10 & 4.976E-10 & 4.758E-10 & 4.581E-10 & 4.449E-10  \\ 
			set14 & 7.708E-10 & 6.491E-10 & 5.772E-10 & 5.436E-10 & 5.179E-10 & 4.978E-10 & 4.898E-10 & 4.627E-10 & 4.470E-10 & 4.442E-10  \\ 
			set15 & 7.883E-10 & 6.848E-10 & 6.463E-10 & 5.897E-10 & 5.469E-10 & 5.256E-10 & 5.274E-10 & 5.102E-10 & 5.084E-10 & 4.854E-10  \\ 
			set16 & 7.417E-10 & 6.373E-10 & 5.863E-10 & 5.656E-10 & 5.243E-10 & 5.344E-10 & 5.094E-10 & 5.068E-10 & 4.806E-10 & 4.738E-10  \\ 
			
			\bottomrule
		\end{tabular*}
	\end{table*}
	
	\begin{table*}[width=2\linewidth,cols=9,pos=h]
		\caption{The dose rate (Sv/t) and spatial utilization rate (\%) under varies running times.}\label{fortime}
		\begin{tabular*}{\tblwidth}{@{} LCCCCCCCCC@{} }
			\toprule
			\multirow{2}{*}{Dataset} & \multicolumn{8}{c}{$\alpha = 0.6$} \\
			\cmidrule{2-9}
			& \multicolumn{2}{c}{30s} & \multicolumn{2}{c}{60s} & \multicolumn{2}{c}{90s} & \multicolumn{2}{c}{120s} \\
			
			\midrule
			set1 & 84.96 & 3.369E-10 & 85.18 & 3.354E-10 & 85.37 & 3.342E-10 & 85.40 & 3.341E-10  \\ 
			set2 & 85.40 & 2.425E-10 & 85.86 & 2.392E-10 & 86.31 & 2.372E-10 & 86.36 & 2.370E-10  \\ 
			set3 & 84.91 & 2.396E-10 & 85.32 & 2.340E-10 & 85.72 & 2.279E-10 & 85.74 & 2.259E-10  \\ 
			set4 & 85.62 & 2.644E-10 & 85.90 & 2.570E-10 & 86.16 & 2.521E-10 & 86.20 & 2.486E-10  \\ 
			set5 & 85.39 & 2.686E-10 & 85.85 & 2.624E-10 & 86.16 & 2.580E-10 & 86.16 & 2.578E-10  \\ 
			set6 & 85.58 & 2.996E-10 & 85.82 & 2.821E-10 & 86.00 & 2.723E-10 & 86.07 & 2.713E-10  \\ 
			set7 & 86.32 & 2.970E-10 & 86.59 & 2.875E-10 & 86.99 & 2.818E-10 & 87.03 & 2.815E-10  \\ 
			set8 & 86.43 & 3.141E-10 & 86.77 & 3.079E-10 & 86.95 & 3.034E-10 & 87.00 & 3.026E-10  \\ 
			set9 & 87.24 & 3.656E-10 & 87.48 & 3.562E-10 & 87.73 & 3.502E-10 & 87.74 & 3.488E-10  \\ 
			set10 & 88.07 & 3.827E-10 & 88.21 & 3.757E-10 & 88.37 & 3.691E-10 & 88.44 & 3.689E-10  \\ 
			set11 & 88.47 & 4.191E-10 & 88.74 & 4.092E-10 & 89.06 & 3.988E-10 & 89.09 & 3.944E-10  \\ 
			set12 & 88.55 & 4.141E-10 & 88.76 & 4.073E-10 & 88.89 & 3.990E-10 & 88.90 & 3.970E-10  \\ 
			set13 & 89.28 & 4.640E-10 & 89.55 & 4.511E-10 & 89.72 & 4.428E-10 & 89.77 & 4.410E-10  \\ 
			set14 & 89.30 & 4.446E-10 & 89.47 & 4.337E-10 & 89.65 & 4.295E-10 & 89.68 & 4.286E-10  \\ 
			set15 & 89.93 & 4.727E-10 & 90.10 & 4.650E-10 & 90.22 & 4.612E-10 & 90.25 & 4.586E-10  \\ 
			set16 & 90.43 & 4.833E-10 & 90.70 & 4.715E-10 & 90.83 & 4.601E-10 & 90.87 & 4.596E-10  \\ 
			
			\bottomrule
		\end{tabular*}
	\end{table*}

	\section{Conclusion}
This paper presents a beam search–based heuristic algorithm for optimizing the loading of nuclear waste boxes into a disposal pool. By carefully designing placement rules, the proposed method generates loading schemes that enhance spatial utilization while minimizing radiation exposure. Additionally, a dataset comprising 1,600 nuclear waste container loading problems is provided to support future research. To our knowledge, this study is the first to address this specific problem. To facilitate further research, the code and dataset used in this study are available at: \url{https://github.com/Yzhjdj/BSNA}.
	
	\bibliographystyle{cas-model2-names.bst}
	
	\bibliography{nuclearBib.bib}
	
\end{document}